\documentclass[11pt]{article}%
\usepackage{amssymb}
\usepackage{amsmath}
\usepackage{amsfonts}
\usepackage{graphicx}
\usepackage{color}
\usepackage[top=2cm,bottom=2cm,left=2.5cm,right=2.5cm]{geometry}%
\providecommand{\U}[1]{\protect\rule{.1in}{.1in}}
\def\R{\mathbb R}
\def\rN{{{\R}^N}}

\def\va{\rightarrow}
\newtheorem{theorem}{Theorem}[section]

\newtheorem{proposition}{Proposition}[section]

\newtheorem{corollary}[theorem]{Corollary}

\newtheorem{remark}[theorem]{Remark}

\numberwithin{equation}{section}
\begin{document}

\title{Estimates for fully anisotropic elliptic equations\\ with a zero order term }
\author{A. Alberico\thanks{Istituto per le Applicazioni del Calcolo \textquotedblleft
M. Picone\textquotedblright (I.A.C.), Sede di Napoli, Consiglio
Nazionale delle Ricerche (C.N.R.), Via P. Castellino 111, 80131
Napoli, Italy. E--mail: a.alberico@iac.cnr.it} -- G. di
Blasio\thanks{Dipartimento di Matematica e Fisica, Universit\`{a}
della Campania Luigi Vanvitelli, Viale Lincon, 5 - 81100 Caserta,
Italy. E--mail: giuseppina.diblasio@unicampania.it} -- F.
Feo\thanks{Dipartimento di Ingegneria, Universit\`{a} degli Studi di
Napoli \textquotedblleft Pathenope\textquotedblright, Centro
Direzionale Isola C4 80143 Napoli, Italy. E--mail:
filomena.feo@uniparthenope.it}}
\date{}\maketitle

\begin{abstract}
Integral estimates for weak solutions to a class of Dirichlet
problems for nonlinear, fully anisotropic, elliptic equations with a
zero order term are obtained using symmetrization techniques.

\end{abstract}

\bigskip


\section{\bigskip Introduction}

We deal with anisotropic Dirichlet problems of the form
\begin{equation}
\left\{
\begin{array}
[c]{ll}%
-\operatorname{div}\left(  a\left(  x,u,\nabla u\right)  \right)
+g(x,u)=f(x) & \mbox{ in }\Omega\\
& \\
u=0 & \mbox{ on }\partial\Omega,
\end{array}
\right.  \label{Prob}%
\end{equation}
where $\Omega$ is a bounded open subset in $\mathbb{R}^{N}$, with
$N\geq2$, 
$a:\Omega\times\mathbb{R}\times\mathbb{R}^{N}\rightarrow\mathbb{R}^{N}$
is a
Carath\'{e}odory function such that, for \emph{a.e.} $x\in\Omega$,%
\begin{itemize}
\item[(H1)] $a(x,\eta,\xi)\cdot\xi\geq\Phi\left(  \xi\right) \qquad \text{ for }
 \left(  \eta,\xi\right)  \in\mathbb{R}\times\mathbb{R}^{N},$  
\end{itemize}
where $\Phi :\mathbb{R}^{N}\rightarrow\left[  0,+\infty\right[ $ is
an $N-$dimensional Young function that does not necessarily verify
the $\Delta_2$-condition (see definition in $\S$2.1 below).
\\
Moreover, we suppose that
$g:(x,s)\in\Omega\times\mathbb{R}\rightarrow \mathbb{R}$ is a
Carath\'{e}odory function satisfying the following conditions:
\begin{itemize}
\item[(H2)]{ $g$ is  a strictly increasing function in $u$ for fixed $x$;}

\item[(H3)]{
for every $r\in \R$,  there exists $h_r \in L^1(\Omega)$ such that
$$|g(x,s)|\leq h_r(x) \qquad \text{ for \emph{a.e.}}\;\; x\in \Omega
\text{  and }\forall s \in \R \text{ with } |s|<r;$$}

\item[(H4)]{
$\qquad\qquad\qquad\quad
g\left(  x,s\right)  s\geq b\left(  s\right)  s \qquad\text{for \emph{a.e.} }%
x\in\Omega \text{  and } \forall s\in\mathbb{R}, $}
\end{itemize}
 where $b:\mathbb{R}\rightarrow\mathbb{R}\quad$is a continuous and
strictly increasing function such that $b\left(  0\right)=0.$
Finally, $f:\Omega \rightarrow\mathbb{R}$ is a nonnegative
measurable function with suitable summability (see $\S$2.3 below).

We stress that the anisotropy of problem \eqref{Prob} is governed by
a general $N-$dimensional convex function of the gradient not
necessarily of polynomial type. In the last three decades, problems
related to differential operators whose growth with respect to the
partial derivatives of $u$ is governed by different powers, are
widely studied (see \emph{e.g.} \cite{AdBF4, AdBF2,
antontsev-chipot-08, BMS, DiNardo-Feo-Guibe,DiNardo-Feo, FGK, FGL,
FS, Gi, Mar}). This interest has led to an extensive investigation
also for problems governed by fully
anisotropic growth conditions (see \cite{A, AC, AdBF1, AdBF3,cianchi anisotropo}%
) and problems related to different type of anisotropy (see
\emph{e.g.} \cite{AFTL, BFK, DPdB}).



Our aim is to prove an estimate in rearrangement form for weak
solutions to problem \eqref{Prob} via symmetrization methods in the
spirit of \cite {cianchi anisotropo}. The symmetrization techniques
enable us to transform the anisotropic problem \eqref{Prob} into the
following isotropic radial problem
\begin{equation}
\left\{
\begin{array}
[c]{ll}%
-\operatorname{div}\left(  \!\dfrac{\Phi_{\blacklozenge}\left(
\left\vert \nabla v\right\vert \right)  }{\left\vert \nabla
v\right\vert ^{2}}\nabla
v\!\right)  +b(v)=f^{\bigstar}(x) & \hbox{in $\Omega^\bigstar$}\\
& \\
v=0 & \hbox{on $\partial\Omega^\bigstar$}
\end{array}
\right.  \label{P sym}%
\end{equation}
defined in the ball $\Omega^{\bigstar}$ centered at the origin and
having the same measure as $\Omega$. Here, $\Phi_{\blacklozenge}$
denotes the Klimov symmetrization of $\Phi$ and $f^{\bigstar}$ is
the symmetric decreasing rearrangement of $f$ (see definition in
$\S$2.2 below). In contrast with the isotropic case (see
\textit{e.g.} \cite{ALT,Diaz breve,Ta1, Ta2,V} and the bibliography
starting with them), in problem \eqref{P sym} not only the domain
and the data are symmetrized but also the ellipticity condition is
subject to an appropriate symmetrization.

In \cite{cianchi anisotropo} it is showed that the symmetric
rearrangement of a solution  $u$ to the anisotropic problem
(\ref{Prob}) is pointwise dominated by the radial solution to the
isotropic problem \eqref{P sym} with no zero order term (\emph{i.e.}
with $b\equiv 0$). Here, we obtain an estimate preserving the
influence of  zero order term. Indeed, we prove an integral (and not
pointwise) estimate between a solution $u$  to problem \eqref{Prob}
and the solution $v$ to problem \eqref{P sym}, \textit{i.e.}
$$\int_\Omega A(b(u))\, dx \leq \int_{\Omega^\bigstar} A(b(v))\, dx$$
  for every one-dimensional Young function $A$.

\medskip
The paper is organized as follows. In Section 2, we recall some
background on Young functions, rearrangements and Orlicz-Sobolev
spaces. In Section 3, we prove the main result and we give some
generalization. Finally, in the Appendix we analyze some questions
related to the existence and the uniqueness of a solution to problem
\eqref{P sym}.

  \section{Preliminaries}

\subsection{Young functions}

Let $\Phi:\mathbb{R}^{N}\rightarrow\left[  0,+\infty\right[  $ be a
(finite-valued) $N-$\textit{dimensional Young function}, namely an
even convex function such that
\[
\Phi\left(  0\right)  =0 \qquad\text{and }\qquad\underset{\left\vert
\xi\right\vert \rightarrow+\infty}{\lim}\Phi\left(  \xi\right)
\text{ }=+\infty.
\]
A standard case considered in literature is when
\begin{equation}
\Phi\left(  \xi\right)  =\underset{i=1}{\overset{N}{{\displaystyle\sum}}%
}\lambda_{i}\left\vert \xi_{i}\right\vert ^{p_{i}}\qquad \text{for
}\xi
\in\mathbb{R}^{N}, \label{Fi pi}%
\end{equation}
for some $\lambda_{i}>0$ and $p_{i}>1$, for any $i=1,\ldots,N.$ An
extension of $\left(  \ref{Fi pi}\right)  $ is given by
\begin{equation}
\Phi\left(  \xi\right)  =\underset{i=1}{\overset{N}{{\displaystyle\sum}}%
}\Upsilon_{i}(\xi_{i})\qquad\text{for }\xi\in\mathbb{R}^{N}, \label{ex2}%
\end{equation}
where $\Upsilon_{i}\,$, for $i=1,\ldots,N$, are one-dimensional
Young functions vanishing only at zero. For instance, we can choose,
for $i=1,\ldots,N$,
\begin{equation}
\Upsilon_{i}(s)=\left\vert s\right\vert ^{p_{i}}\left(  \log\left(
c+\left\vert s\right\vert \right)  \right)
^{\alpha_{i}}\qquad\text{for } s\in\mathbb{R}\, , \label{ex3}
\end{equation}
 where either $p_{i}>1$ and $\alpha_{i}\in\mathbb{R}$ or $p_{i}=1$
and $\alpha_{i}\geq0$, and the constant $c$ is large enough for
$\Upsilon_{i}$ to be convex.

If we define the \textit{Young conjugate}  of $\Phi$ the function
$\Phi_{\bullet}\left(  \xi^{\prime}\right)  =\sup\left\{
\xi\cdot\xi^{\prime }-\Phi\left(  \xi\right)
:\xi\in\mathbb{R}^{N}\right\} \hbox{ for}\;\;
\xi^{\prime}\in\mathbb{R}^{N},
$
the \textit{Young inequality} tells us that
\begin{equation}
\xi\cdot\xi^{\prime}\leq\Phi\left(  \xi\right) +\Phi_{\bullet}\left(
\xi^{\prime}\right)\qquad \text{for \
}\xi,\xi^{\prime}\in\mathbb{R}^{N},
\label{Young}%
\end{equation}
where $\lq \lq \,\cdot\,"$  stands for scalar product in $\rN$.
We observe that the function $\Phi_{\bullet}$ enjoys the same
properties as $\Phi$ and is a (finite-valued) $N-$dimensional Young
function if
\begin{equation}
\underset{\left\vert \xi\right\vert
\rightarrow+\infty}{\lim}\frac{\Phi\left(
\xi\right)  }{\left\vert \xi\right\vert }=+\infty. \label{Young-func_2}%
\end{equation}

 The following result (see \cite[Theorem
5.1]{S1}) says us when the Young inequality \eqref{Young} holds as
equality.
\begin{proposition}
\label{equalityYoung} If the derivative of $\Phi$
with respect to a direction $\eta$ is linear in $\eta,$ then, given
any $\xi_{0}$, the i-th component of $\eta$, $\eta_{i}$, verifies
$\eta_{i}=\frac{\partial \Phi}{\partial e_{i}}(\xi _{0})$ if and
only if $\xi_{0}\cdot\eta=\Phi\left( \xi_{0}\right)  \Phi_{\bullet
}\left( \eta\right)$. Here, $\left\{ e_{i}\right\}$ stands for a
basis in ${{\mathbb{R}}^{N}.}$
\end{proposition}
Finally, we recall a useful property of an $N$-dimensional Young
function that we will use in what follows  (see \cite[Proposition
6.7]{CB}).
\begin{proposition}\label{Young_equ}
Let $\Phi \in C^1 (\R^N)$ be an $N$-dimensional Young function such
that \eqref{Young-func_2} holds. Then,
\begin{equation}\label{dis CB}
\Phi_{\bullet }(\Phi_{\xi }(\xi))\leq \Phi_{\xi }(\xi) \cdot \xi
\leq \Phi(2\xi)\, .
\end{equation}
\end{proposition}
An $N-$dimensional Young function $\Phi$ is said to satisfy the
\emph{$\Delta_2-$condition} near infinity if  there exist constants
$C>1$ and $K\geq0$ such that $\Phi(2\xi)\leq C\, \Phi(\xi)$ for
$|\xi|>K$. For example, one can easily verify that  in \eqref{ex2}
every function $\Upsilon_i$ satisfies the $\Delta_2-$condition
whenever $\Phi(\xi)$ does. Conversely, if one function $\Upsilon_i$
does not satisfy the $\Delta_2-$condition (as it happens for
$\Upsilon_i$ given by \eqref{ex3}), then $\Phi(\xi)$ does the same.
Another example of a function which does not satisfy the
$\Delta_2-$condition is given by
$\Phi(\xi)=e^{\sum_{i=1}^N|\xi_i|^{p_i}}-1$ with $p_i\geq1$ for any
$i=1,\cdots,N$. Indeed,
$$\lim_{|\xi|\rightarrow
+\infty}\Big(\sum_{i=1}^N|2\xi_i|^{p_i}-\sum_{i=1}^N|\xi_i|^{p_i}\Big)=
\lim_{|\xi|\rightarrow
+\infty}\sum_{i=1}^N|\xi_i|^{p_i}(2^{pi}-1)=+\infty.$$

\subsection{Symmetrization}

A precise statement of our result requires the use of classical
notions of rearrangement of a function (see \textit{e.g.} \cite{BS})
and of suitable symmetrization of a Young function introduced by
Klimov in \cite{Klimov 74}. \newline Let $u$ be a measurable
function (continued by $0$ outside its domain) fulfilling
\begin{equation}
\left\vert \{x\in\mathbb{R}^{N}:\left\vert u(x)\right\vert
>t\}\right\vert <+\infty \qquad\text{ for every } t>0.
\label{insieme livello di misura finita}%
\end{equation}
The \textit{symmetric decreasing rearrangement} of $u$ is the
function $u^{\bigstar}:\mathbb{R}^{N}\rightarrow\left[
0,+\infty\right[  $ $\ $satisfying
\begin{equation}
\{x\in\mathbb{R}^{N}:u^{\bigstar}(x)>t\}=\{x\in\mathbb{R}^{N}:\left\vert
u(x)\right\vert >t\}^{\bigstar}\qquad\text{for }t>0.\label{livello palla}%
\end{equation}
The \textit{decreasing rearrangement} $u^{\ast}$ of $u$ is defined
as
\[
u^{\ast}(s)=\sup\{t>0:\mu_{u}(t)>s\}\qquad\text{for }s\geq0,
\]
where $ \mu_{u}(t)=\left\vert \{x\in{\Omega}:\left\vert
u(x)\right\vert
>t\}\right\vert \text{ for }t\geq0
$ denotes the \textit{distribution function} of $u$. Moreover,
\[
u^{\bigstar}(x)=u^{\ast}(\omega_{N}\left\vert x\right\vert ^{N})
\qquad\text{ for \emph{a.e.}}\;x\in{{\mathbb{R}}^{N},}%
\]
where $\omega_{N}$ is the measure of the $N-$dimensional unit ball.
Analogously, we define the \textit{symmetric increasing
rearrangement} $u_{\bigstar}$ of $u$ on replacing
\textquotedblleft$>$\textquotedblright\ by
\textquotedblleft$<$\textquotedblright\ in the definitions of the
sets in (\ref{insieme livello di misura finita}) and (\ref{livello
palla}).


Let $\Phi$ be an $N-$dimensional Young function. We denote by $\Phi
_{\blacklozenge}:\mathbb{R\rightarrow}\left[  0,+\infty\right[  $
the symmetrization of $\Phi$ introduced in \cite{Klimov 74} (see
also \cite{CianchiFully, cianchi anisotropo}). It is the
one-dimensional Young function defined by
\begin{equation*}
\Phi_{\blacklozenge}(\left\vert \xi\right\vert
)=\Phi_{\bullet\bigstar\bullet
}\left(  \xi\right) \qquad\text{for }\xi\in\mathbb{R}^{N}, \label{def fi rombo}%
\end{equation*}
namely it is the composition of Young conjugation, symmetric
increasing rearrangement and Young conjugate again. We stress that
the functions $\Phi_{\blacklozenge}$ and $\Phi_{\bigstar}$ are not
equals in general, but they are always equivalent, \emph{i.e.} there
exist two positive constants $K_{1}$ and $K_{2}$ such that
\begin{equation}
\Phi_{\bigstar}(K_{1}\xi)\leq\Phi_{\blacklozenge}(|\xi|)\leq\Phi_{\bigstar
}(K_{2}\xi)\qquad\hbox{\rm for}\;\;\xi\in{{\mathbb{R}}^{N}}. \label{equiv}%
\end{equation}
Nevertheless, $\Phi_{\blacklozenge}(\left\vert \,\cdot\,\right\vert
)=\Phi _{\bigstar}\left(  \,\cdot\,\right)  $ if and only if $\Phi$
is radial,\textit{ i.e.} $\Phi=\Phi_{\bigstar}.$

\noindent We emphasize that, if $\Phi$  verifies
\eqref{Young-func_2}, then
\begin{equation}\label{101}
\underset{s \rightarrow +\infty}{\lim}\frac{\Phi_\blacklozenge(s)}
{s}=+\infty.
\end{equation}

 \noindent Moreover, the Young function
$\Phi_{\blacklozenge}$ verifies the $\triangle_2$-condition whenever
$\Phi$ does. For the convenience of the reader we give the details.
Since $\Phi(\xi)\rightarrow +\infty$ as $|\xi|\rightarrow + \infty$,
there exists $M>0$ such that $\{\xi : \Phi(\xi)> t\}\cap\{\xi :
|\xi|>M\}=\{\xi : \Phi(\xi)> t\}$ for every $ t>0$. The function
$\Phi(\xi)$ verifies the $\triangle_2$-condition if and only if
$\{\xi : \Phi(\xi)> \frac{t}{c}\}\subseteq \{\xi : \Phi(2\xi)> t\}$.
The assert follows observing the equimisurability of $\Phi$ and
$\Phi_{\bigstar}$, the equality $\{\xi :
\Phi_{\bigstar}(\xi)>\tau\}=\mathbb{R}^N \setminus B(0,R)$ with
$R^N=\frac 1{\omega_N}|\{\xi : \Phi(\xi)> \tau\}|$ and the
equivalence between $\Phi_{\blacklozenge}$ and $\Phi_{\bigstar}$.
\newline When $\Phi$ is  given by \eqref{Fi pi}, easy calculations
show that
\begin{equation*}
\Phi_{\blacklozenge}(\left\vert \xi\right\vert )=\Lambda\left\vert
\xi\right\vert ^{\overline{p}}, \label{fi rombo}%
\end{equation*}
where $\overline{p}$ is the harmonic mean of exponents
$p_{1},\ldots,p_{N}$
 and $\Lambda$ is a suitable positive constant.
\newline In the more general case \eqref{ex2}, it is possible to
show that
\begin{equation}
\Phi_{\blacklozenge}^{-1}(r)\approx\left(  \underset{i=1}{\overset{N}{\Pi}%
}\Upsilon_{i}^{-1}(\cdot)\right)  ^{\frac{1}{N}}(r)\qquad \text{ for
}r\geq0,
\label{inversa fi rombo}%
\end{equation}
where $\Upsilon_{i}^{-1}$ denotes the inverse function of
$\Upsilon_{i}$ in $\left[  0,+\infty\right)$.
From (\ref{inversa fi rombo}), we deduce that%
\begin{equation*}
\Phi_{\blacklozenge}(s)\approx\left\vert s\right\vert
^{\overline{p}}\left(
\log\left(  c+\left\vert s\right\vert \right)  \right)  ^{\frac{\overline{p}%
}{N}\underset{i=1}{\overset{N}{%
{\textstyle\sum}
}}\frac{\alpha_{i}}{p_{i}}} \label{fi rombo con log}%
\end{equation*}
near infinity when $\Upsilon_{i}$ is given by \eqref{ex3}.
\newline

\noindent
 Now, let us define  the following function
\begin{equation}
\Psi_{\blacklozenge}(s)=\left\{
\begin{array}
[c]{ll}%
\displaystyle\frac{\Phi_{\blacklozenge}(s)}{s} & \hbox{\rm for $s >0 $}\\
& \\
0 & \hbox{\rm for $s=0$.}
\end{array}
\right.  \label{Psi}%
\end{equation}
It is not-decreasing and, if, in addiction,%
\begin{equation}
\underset{s\rightarrow0^{+}}{\lim}\frac{\Phi_{\blacklozenge}\left(
s\right)
}{s}=0, \label{psi incrising}%
\end{equation}
then it is strictly increasing in $\left[ s_{0},+\infty\right)  $
with
\begin{equation}
s_{0}=\sup\left\{  s\geq0:\Phi_{\blacklozenge}\left(  s\right)
=0\right\}. \label{s0}
\end{equation}
 Moreover,
\begin{equation}
\Phi_{\blacklozenge\bullet}\left(  r\right)
\leq\Phi_{\blacklozenge}\left( \Psi_{\blacklozenge}^{-1}\left(
r\right)  \right) \qquad \text{for \ \ }r\geq0,
\label{dis PsiRombo}%
\end{equation}
where $\Psi_{\blacklozenge}^{-1}$ denotes the inverse of
$\Psi_{\blacklozenge }$ restricted to $[s_{0},+\infty)$.

We recall that in the anisotropic setting a \textit{Polya-Szeg\"{o}
principle} holds (see \cite{cianchi anisotropo}). More precisely,
let $u$ be a weakly differentiable function
in $%
\mathbb{R}
^{N}$ satisfying (\ref{insieme livello di misura finita}) and such
that $\int_{\mathbb{R}^{N}}\Phi\left(  \nabla u\right)  dx<+\infty$.
Then,
$u^{\bigstar}$ is weakly differentiable in $%
\mathbb{R}
^{N}$ and
\begin{equation}
\int_{\mathbb{R}^{N}}\Phi_{\blacklozenge}\big(  \big\vert \nabla
u^{\bigstar}\big\vert \big)\; dx\leq\int_{\mathbb{R}^{N}}\Phi\left(
\nabla
u\right)\;  dx . \label{cianchi}%
\end{equation}


\subsection{\label{Spases}Weak solutions and Orlicz Spaces}
In according to the ellipticity condition (H1), we look for weak
solution $u$ to problem \eqref{Prob} in the class
$V_{0}^{1,\Phi}(\Omega)$ of real-valued functions $u$ in $\Omega$,
whose continuation by $0$ outside $\Omega$ is weakly differentiable
in $\mathbb{R}^{N}$ and satisfies $\int_{\Omega} \Phi\left(  \nabla
u\right)  dx<\infty $.
Let us suppose that the datum $f$ is such that
\begin{equation}\label{f}\int_\Omega\Phi_{\bullet}(f(x))\,
dx<+\infty.\end{equation} A function $u\in V_{0}^{1,\Phi}(\Omega)$
is a weak solution to problem (\ref{Prob}) if $b(u)\in
L^{1}(\Omega),u\,b(u)\in L^{1}(\Omega)$ and
\begin{equation}
\int_{\Omega}a\left(  x,u,\nabla u\right)  \cdot\nabla\varphi
    \; dx+\int_{\Omega
}g(x,u)\,\varphi\;dx=\int_{\Omega}f\,\varphi\;dx\text{ }\ \  \label{sol deb}%
\end{equation}
for every $\varphi\in V_{0}^{1,\Phi}(\Omega)\cap L^{\infty}(\Omega)$
and for $\varphi=u$.

\noindent We recall that\ $V_{0}^{1,\Phi}(\Omega)$ is always a
convex set, but not necessary a linear space since no
$\Delta_{2}-$condition on $\Phi$ is required.

\medskip
Let $\Omega$ be a measurable subset in ${{\mathbb{R}}^{N}}$ having
finite Lebesgue measure, we say that  $A:[0, \infty)\va [0, \infty)$
is a one$-$dimensional Young function if it is a convex,
left-continuous function, vanishing at $0$, and neither identically
equal to $0$, nor to $+\infty$.  The \emph{Orlicz space}
$L_{A}(\Omega)$, associated with the Young function $A$, is the set
 of all measurable functions $g$ in
$\Omega$ for which the Luxemburg norm
\[
\| g\|_{L_{A}(\Omega)}=\inf\left\{ \lambda>0:\int_{\Omega}A\left(
\frac{|g(x)|}{\lambda}\right) dx\leq1\right\}
\]
is finite. The functional $\|\cdot \|_{L_{A}(\Omega)}$ is a norm on
$L_A(\Omega)$ which renders the latter a Banach function space.
A generalized H\"{o}lder inequality
\[
{\displaystyle\int_{\Omega}}|f\, g|\; dx\leq2\left\Vert f\right\Vert
_{L_{A}(\Omega)}\left\Vert g\right\Vert _{L_{A_{\bullet}}(\Omega)}\text{ }%
\]
holds for every $f\in L_{A}(\Omega)$ and $g\in
L_{A_{\bullet}}(\Omega)$, where $A_{\bullet}$ is the Young conjugate
of $A$.
\\
We define the space $E_{A}(\Omega)$ as the closure of
$L^{\infty}(\Omega )\;$in$\;L_{A}(\Omega).$ In general
$E_{A}(\Omega)\subset L_{A}(\Omega)$ unless $A$ satisfies the
$\Delta_{2}-$condition. The space $L_{A}(\Omega)$ is the dual space
of $E_{A_{\bullet}}(\Omega)$ and
the duality pairing is given by%
\[
<f,g>={\displaystyle\int_{\Omega}}f\;g\;dx
\]
for $f\in L_{A}(\Omega)$ and $g\in E_{A_{\bullet}}(\Omega)$. 
\\
 The \emph{Orlicz-Sobolev} space $W^{1}L_{A}(\Omega)$  is defined
requiring that $u$ and $\left\vert \nabla u\right\vert$ are in
$L_{A}(\Omega)$.
It is a Banach function space equipped with the norm $\left\Vert u\right\Vert _{W^{1}%
L_{A}(\Omega)}=\left\Vert |\nabla u| \right\Vert _{L_{A}(\Omega)}.$
\\
 Finally, we define $W_{0}^{1}L_{A}(\Omega)$ as the closure of
$\mathcal{D}(\Omega)$ in $W^{1}L_{A}(\Omega)$ with respect to the
weak topology $\sigma(L_{A}(\Omega),E_{A_{\bullet}}(\Omega)).$ For
more details about Orlicz spaces, we refer to \cite{rao}.

Now we recall the \emph{embedding theorem for anisotropic
Orlicz-Sobolev spaces} (see \cite{CianchiFully}). Assume that $\Phi$
is an $N-$dimensional Young function fulfilling
\begin{equation*}
\int_{0}\left(  \frac{s}{\Phi_{\blacklozenge}(s)}\right)  ^{\frac{1}{N-1}%
}\;ds<\infty\, ,\label{cond_in zero}%
\end{equation*}
and let $\Phi _{N}:[0,\infty)\rightarrow\lbrack0,\infty]$ be the
\emph{Sobolev conjugate} of $\Phi$ defined as
\begin{equation*}
\Phi_{N}(s)=\Phi_{\blacklozenge}(H_{\Phi_{\blacklozenge}}^{-1}(s))\quad
\hbox{for}\;\; s\geq 0 \quad \text{ where } \quad
H_{\Phi_{\blacklozenge}}(r)=\left( \int_{0}^{r}\left(
\frac{s}{\Phi_{\blacklozenge}(s)}\right) ^{\frac{1}{N-1}}ds\right)
^{\frac{1}{N^{\prime}}} \quad \hbox{\rm for}\;\; r\geq 0\, .
\end{equation*}
\\
 If
\begin{equation}
\int^{+\infty}\left(  \frac{s}{\Phi_{\blacklozenge}(s)}\right)
^{\frac
{1}{N-1}} \; ds<\infty\, , \label{existence_1}%
\end{equation}
then any function $u\in V_{0}^{1,\Phi}(\Omega)$ is essentially
bounded.
\\
If, instead
\begin{equation}
\int^{+\infty}\left(  \frac{s}{\Phi_{\blacklozenge}(s)}\right)
^{\frac
{1}{N-1}}\; ds=\infty \, , \label{existence_2}%
\end{equation}
then
\[
\int_{\Omega}\Phi_{N}(c\,u(x))dx<\infty
\]
for any function $u\in V_{0}^{1,\Phi}(\Omega)$ and for every $c\in
\R$.

We stress that we are not interested in the existence of solution to
problem (\ref{Prob}). Nevertheless, thanks to the embedding theorem
for anisotropic Orlicz-Sobolev spaces above, we give some conditions
on $a(x,\eta,\xi)$ and $f\left(  x\right)  $ in order to guarantee
that (\ref{sol deb}) is well posed. More precisely, assume that
condition \eqref{existence_1} holds. Then, the left-hand side of
(\ref{sol deb}) is always finite and the right-hand side of
(\ref{sol deb}) is finite if $f\in$ $L^{1}\left( \Omega\right)$.
Otherwise, if condition \eqref{existence_2} holds, we shall require
some conditions on functions $a$ and $f$.
 More precisely,
\begin{equation}\label{ff}
\|s^{\frac{1}{N}}f^{\ast\ast}\left(  s\right)\|_{L_{\Phi_{\blacklozenge\bullet}%
}\left(  0,\left\vert \Omega\right\vert \right)}<\infty
\end{equation}
and
\[
\Phi_{\bullet}(a(x,\eta,\xi))\leq c\left[  \theta\left(  x\right)
+M\left( \eta\right)  +\Phi(\xi)\right]
\qquad\text{for}\;\;a.e.\;\;x\in \Omega\;  \;\text{and for}\;\left(
\eta,\xi\right)  \in\mathbb{R}\times \mathbb{R}^{N},
\]
where $ f^{\ast\ast}(s)=\frac{1}{s}\int_{0}^{s}f^{\ast}(r)\;dr $,
$c$ is a positive constant, $\theta$ is a positive function in
$L^{1}\left(  \Omega\right)$   and $M:\mathbb{R\rightarrow}\left[
0,+\infty\right[  $ is a continuous function such that $M\left(
\eta\right) \leq\Phi_{N}\left( k\eta\right)  $ for some $k>0$ and
for every $\eta\in\mathbb{R}.$ Note that condition  \eqref{ff} is
weaker than  \eqref{f}.


\section{Mass comparison results}

\subsection{\bigskip Symmetrized problem}

Our aim is to prove a comparison between the concentration of a
solution $u$ to problem \eqref{Prob} and the solution $v$ to a
suitable simpler isotropic problem. More precisely, in our case the
so-called ``symmetrized'' problem is
\begin{equation}
\left\{
\begin{array}
[c]{ll}%
-\operatorname{div}\left(  \!\dfrac{\Phi_{\blacklozenge}\left(
\left\vert \nabla v\right\vert \right)  }{\left\vert \nabla
v\right\vert ^{2}}\nabla
v\!\right)  +b(v)=\widetilde{f}(x) & \hbox{in $\Omega^\bigstar$}\\
& \\
v=0 & \hbox{on $\partial\Omega^\bigstar$}\,,
\end{array}
\right.  \label{Prob_sym_zero}%
\end{equation}
where

\begin{itemize}
\item[(H$^\bigstar1$)]  $ \Phi_{\blacklozenge}(t)$ is a one-dimensional
 Young function strictly increasing in $[0,+\infty)$ such that
 \eqref{101}
 and \eqref{psi incrising} hold,

    \item[(H$^\bigstar2$)] $b$ is a continuous strictly increasing function such that $b(0)=0$,

 \item[(H$^\bigstar3$)] $\Omega^{\bigstar}$ is the ball centered at the origin having the same
measure as $\Omega$,

\item[(H$^\bigstar4$)] $ \,\widetilde{f}:\Omega^{\bigstar}\rightarrow\mathbb{R}$ is a
nonnegative radially symmetric and decreasing along the radii,

\item[(H$^\bigstar5$)] $  \widetilde{f}\in E_{\Theta_{\bullet}}\left(  \Omega^{\bigstar}\right)
,$ where $\Theta:\mathbb{R}\rightarrow\lbrack0,\infty)$ is the Young
function defined by
\begin{equation}
\Theta\left(  r\right)
=\int_{0}^{|r|}\frac{\Phi_{\blacklozenge}\left(
s\right)  }{s}ds\qquad\hbox{for $r \in \mathbb{R} $}\,. \label{theta}%
\end{equation}

\end{itemize}

By Proposition 13 of \cite{GM}, under the previous assumptions,
there exists a unique positive weak solution $v\in
W_{0}^{1}L_{\Theta}(\Omega^{\bigstar})\bigcap\left\{v \in
W_{0}^{1}L_{\Theta}(\Omega^{\bigstar}):
\left\vert\frac{\Phi_{\blacklozenge}(|\nabla v|)}{|\nabla
v|^2}\nabla v \right \vert \in
L_{\Theta_{\bullet}}(\Omega^{\bigstar})\right\} $ to problem
\eqref{Prob_sym_zero} (see also
\cite{Diaz breve} and the Appendix for more details and remarks),
\textit{i.e.}
$$\int_{\Omega^{\bigstar}}\left[\frac{\Phi_{\blacklozenge}(|\nabla v|)}
{|\nabla v|^2}\nabla v \cdot\nabla\varphi+ b(v) \varphi \right]
\,dx= \int_{\Omega^\bigstar}\widetilde{f}\, \varphi \; dx\qquad
\text{for every}\;\; \varphi \in
W_{0}^{1}L_{\Theta}(\Omega^{\bigstar})\cap
L^{\infty}(\Omega^{\bigstar}).$$
 We stress that $\left
\vert\frac{\Phi_{\blacklozenge}(|\nabla v|)}{|\nabla v|^2}\nabla v
\right \vert$ does not necessarily belong to the space
$L_{\Theta_{\bullet}}(\Omega^{\bigstar})$ for every $v\in
W_{0}^{1}L_{\Theta}(\Omega^{\bigstar})$, since $\Theta$ does not
necessarily verify a $\Delta_2$-condition.
 Indeed, the function $\Theta$ verifies the $\Delta_2$-condition whenever $\Phi_\blacklozenge$ does
 \footnote{
 Fixed $t_1\in \mathbb{R}$, using the $\Delta_2$-condition for $\Phi_\blacklozenge$,  we get, for
 $|\xi|>2t_1$,
 \begin{align*}
 \Theta(2|\xi|)=\int_0^{2t_1}\frac{\Phi_\blacklozenge(s)}{s}ds+\int_{t_1}^{|\xi|}\frac{\Phi_\blacklozenge(2s)}{s}ds
 \leq  \Theta(|\xi|)+ K_1 \Theta(|\xi|)
 \end{align*}
 for some $K_1>0$.
 Otherwise, since $\Theta$ is $C^1$ and strictly increasing, fixed $t_2\in \mathbb{R}$, we obtain, for
 $|\xi|>t_2$,
$$
 \frac{\Theta(2|\xi|)}{\Theta(|\xi|)}=\frac{\Theta(2|\xi|)-\Theta(2t_2)}
 {\Theta(|\xi|)-\Theta(t_2)}\frac{\Theta(2|\xi|)}{\Theta(2|\xi|)-\Theta(2t_2)}
 \frac{\Theta(|\xi|)-\Theta(t_2)}{\Theta(|\xi|)} \geq  \frac{\Phi_\blacklozenge(2c)}
 {\Phi_\blacklozenge(c)}(1-\epsilon)
 $$
for some $\epsilon>0$ and $t_2<c<|\xi|$. If $\Delta_2$-condition is
not fulfilled, then $
 \frac{\Theta(2|\xi|)}{\Theta(|\xi|)}>K_2$, for $|\xi|>t_2$  and for every $K_2 \in \mathbb{R}.$
}.

We look for solutions in the class $
V_{0}^{1,\Phi_{\blacklozenge}}\left( \Omega^{\bigstar}\right)$.
Using the monotonicity of $\Psi_\blacklozenge(r)$, defined as in
\eqref{Psi}, we get
\begin{equation}\label{theta2}
\Theta(r) \leq \Phi_\blacklozenge (|r|) \qquad\hbox{\rm for}\;\;
r\in \R.
\end{equation}
For this reason we need to require more summability on datum
$\widetilde{f}$ with to respect to hypothesis (H$^\bigstar$5) in
order to assure that the weak solution $v\in
W_{0}^{1}L_{\Theta}(\Omega^{\bigstar})\bigcap\left\{v \in
W_{0}^{1}L_{\Theta}(\Omega^{\bigstar}):
\left\vert\frac{\Phi_{\blacklozenge}(|\nabla v|)}{|\nabla
v|^2}\nabla v \right \vert \in
L_{\Theta_{\bullet}}(\Omega^{\bigstar})\right\} $ belongs to the
class $V_{0}^{1,\Phi_{\blacklozenge}}\left(
\Omega^{\bigstar}\right).$ More precisely, if we require the
following additional assumptions
\begin{itemize}
\item[(H$^\bigstar$6)]{
$
\text{either}\quad
\underset{r\rightarrow+\infty}{\lim}\Psi_{\blacklozenge}\left(
r\right) =+\infty
\quad \text{ or} \quad
\dfrac{s^{1/N}}{N\omega_{N}^{1/N}}\widetilde{f}^{\ast\ast}(s)
<\text{\ }\underset{r\rightarrow+\infty}{\lim}\Psi_{\blacklozenge
}\left( r\right)  \text{ \ \ \ for every }s>0,
$}
\end{itemize}

\begin{itemize}
 \item[(H$^\bigstar$7)]{
$ \displaystyle\int_{0}^{\left\vert \Omega\right\vert
}\Phi_{\blacklozenge }\left( \Psi_{\blacklozenge}^{-1}\left(
\dfrac{s^{1/N}}{N\omega_{N}^{1/N}}\widetilde
{f}^{\ast\ast}(s)\right) \right) ds\quad \text{  is finite,} $}
\end{itemize}
 then $v \in V_{0}^{1,\Phi_{\blacklozenge}}\left(
\Omega^{\bigstar}\right).$ Indeed, condition (H$^\bigstar$6),
monotonicity of $\Psi
_{\blacklozenge}^{-1}$ and (H$^\bigstar$7) allow to obtain%
\begin{align}
\displaystyle\int_{\Omega^{\bigstar}}\Phi_{\blacklozenge}\left(
\left\vert \nabla v\right\vert \right)  dx &=\int_{0}^{\left\vert
\Omega^{\bigstar }\right\vert }\Phi_{\blacklozenge}\left(
\Psi_{\blacklozenge}^{-1}\left( \dfrac{\int_{0}^{r}\left[
\widetilde{f}^{\ast}(s)-b\left(  v^{\ast
}\left(  s\right)  \right)  \right]  ds}{N\omega_{N}^{1/N}r^{1/N^\prime}%
}\right)  \right)  dr
\cr & \leq\int_{0}^{\left\vert \Omega^{\bigstar}\right\vert
}\Phi_{\blacklozenge}\left(  \Psi_{\blacklozenge}^{-1}\left(
\dfrac{r^{1/N}\widetilde
{f}^{\ast\ast}(r)}{N\omega_{N}^{1/N}}\right)  \right)  dr<+\infty.
\nonumber
\end{align}
Using  (\ref{dis PsiRombo}), it is easy to check that
 condition (H$^\bigstar$7) implies \eqref{ff}
with $f=\widetilde{f}$ and that
$\Phi_{\blacklozenge\bullet}\left(\left
\vert\frac{\Phi_{\blacklozenge}(|\nabla v|)}{|\nabla v|^2}\nabla v
\right \vert\right) $ belongs to $L^1(\Omega^{\bigstar})$, for every
$v\in V_0^{1,\Phi_{\blacklozenge}}(\Omega^\bigstar)$. Moreover,
using also \eqref{theta2}, we have that
$$\int_{\Omega^{\bigstar}}\left[\frac{\Phi_{\blacklozenge}
(|\nabla v|)}  {|\nabla v|^2}\nabla v \cdot \nabla\varphi+ b(v)
 \varphi \right] \,dx= \int_{\Omega^\bigstar}\widetilde{f}\,
 \varphi\; dx
 \qquad \text{holds\;\, for}\;\; \varphi \in  V_{0}^{1,\Phi_{\blacklozenge}}
 \left( \Omega^{\bigstar}\right)\cap L^{\infty}(\Omega^{\bigstar}).$$

Finally, the symmetry of data assures that the solution $v$ of
(\ref{Prob_sym_zero}) is radially symmetric and moreover $v\left(
x\right)  =v^{^{\bigstar}}\left(  x\right)  $ (for more details, see
the proof of Theorem \ref{teomain} below).

\subsection{\bigskip Statements}

Let us consider problems \eqref{Prob} and \eqref{Prob_sym_zero}. In
the following theorem, we prove an $L^{\infty}-$norm estimate of the
difference between the concentration of solutions in terms of
difference between the concentration of data. Moreover, under the
additional assumption that datum $f$ is less concentrated than datum
$\widetilde{f},$ we obtain a comparison between the concentration of
a solution $u$ to problem (\ref{Prob}) and the concentration of the
solution $v$ to problem (\ref{Prob_sym_zero}).

\begin{theorem}
\label{teomain} Let $\Phi:\mathbb{R}^{N}\rightarrow\left[
0,+\infty\right[  $ be an $N-$dimensional Young function fulfilling
\eqref{Young-func_2}, such that $\Phi_{\blacklozenge}(t)$ is
strictly increasing in $[0,+\infty)$ verifying  \eqref{psi
incrising}. Assume that \rm{(H1)-(H4), (H$^\bigstar$3),
(H$^\bigstar$4), (H$^\bigstar$6), (H$^\bigstar$7)} \emph{hold and
$f$   nonnegative. If $u$ is a (nonnegative) weak solution to
\eqref{Prob} and $v$  the (nonnegative) weak solution to
\eqref{Prob_sym_zero}, then
\begin{equation*}
\Vert(\mathcal{B}-\widetilde{\mathcal{B}})_{+}\Vert_{L^{\infty}(0,|\Omega
|)}\leq\Vert(\mathcal{F}-\widetilde{\mathcal{F}})_{+}\Vert_{L^{\infty
}(0,|\Omega|)},\label{3:10}%
\end{equation*}
where $(\tau)_+:=\max\{0, \tau\}$, and
\begin{align}
&  \mathcal{B}(s)=\int_{0}^{s}b(u^{\ast}(t))\;dt & \quad &
\widetilde
{\mathcal{B}}(s)=\int_{0}^{s}b(v^{\ast}(t))\;dt\label{3:11}\\
&  \mathcal{F}(s)=\int_{0}^{s}f^{\ast}(t)\;dt & \quad &  \widetilde
{\mathcal{F}}(s)=\int_{0}^{s}\widetilde{f}^{\ast}(t)\;dt\label{3:11'}%
\end{align}
for $s\in(0,|\Omega|]$.}
\end{theorem}

\medskip

An immediate consequence of Theorem \ref{teomain} \ is the following
integral estimate of $b(u)$ in terms of  $b(v)$.

\begin{corollary}
Under the same assumptions of Theorem \ref{teomain}, if
\begin{equation}\label{100}
\mathcal{F}(s)\leq\widetilde{\mathcal{F}}(s)\qquad\hbox{for
any}\,\;s\in\lbrack0,|\Omega|,
\end{equation}
then%
\[
\mathcal{B}(s)\leq\widetilde{\mathcal{B}}(s)\qquad\hbox{for
any}\,\;s\in\lbrack0,|\Omega|].
\]
Moreover, \begin{equation*} \int_{\Omega}A(b
(u(x)))\;dx\leq\int_{\Omega^{\bigstar}}A
(b(v(x)))\;dx\label{norm estim}%
\end{equation*}
for all convex and non-decreasing function $A:[0,+\infty)\rightarrow
[0,+\infty)$ such that $A(0)=0.$
\\
If $b$ is strictly increasing, then%
\[
\left\Vert u\right\Vert _{L^{\infty}(\Omega)}\leq\left\Vert
v\right\Vert
_{L^{\infty}(\Omega^{\bigstar})}.%
\]

\end{corollary}

\smallskip
\noindent
 Analogous results, for $\Phi$ given by \eqref{Fi pi}, are
contained in \cite{AdBF3, AdBF4}.

\noindent We note that if $\widetilde{f}= f^\bigstar$, then
\eqref{100} holds.

 \begin{remark} \rm The use of Klimov symmetrization
gives sharp estimates in previous results. Indeed, in the proof of
Theorem \ref{teomain}, replacing  $\Phi_\blacklozenge$ by
$\Phi_\bigstar$ and using \eqref{equiv}, it follows that there
exists a constant $K>0$ such that
\begin{equation}\label{K}
\|(K \mathcal{B} -
\widetilde{\mathcal{B}})_+\|_{L^\infty(\Omega)}\leq \| (K
\mathcal{F} - \widetilde{\mathcal{F}})_+\|_{L^\infty(\Omega)}\,.
\end{equation}
In particular, if there exists $h>0$ such that
\begin{equation*}
\mathcal{F}(s)\leq h\, \widetilde{\mathcal{F}}(s)\qquad\hbox{for
any}\,\;s\in\lbrack0,|\Omega|],\label{3:12}%
\end{equation*}
then
\[
\mathcal{B}(s)\leq h \,\widetilde{\mathcal{B}}(s)\qquad\hbox{for
any}\,\;s\in\lbrack0,|\Omega|].
\]

\end{remark}

\subsection{Proof of Theorem \ref{teomain}}


We define the functions $u_{\kappa,t}:\Omega\rightarrow$\textbf{
}$\mathbb{R}$ as
\[
u_{\kappa,t}\left(  x\right)  =\left\{
\begin{array}
[c]{ll}%
0 & \mbox{ if }\left\vert u\left(  x\right)  \right\vert \leq
t,\vspace
{0.2cm}\\
\left(  \left\vert u\left(  x\right)  \right\vert -t\right)  \text{sign}%
\left(  u\left(  x\right)  \right)  & \mbox{ if
}t<\left\vert u\left(  x\right)  \right\vert \leq t+\kappa\\
& \\
\kappa\;\text{sign}\left(  u\left(  x\right)  \right)  & \mbox{ if
}t+\kappa<\left\vert u\left(  x\right)  \right\vert
\end{array}
\right.
\]
for any fixed $t$ and $\kappa>0.$ This function can be chosen as
test function in (\ref{sol deb}). Using \eqref{cianchi} and arguing
as in \cite{cianchi anisotropo}, we have
\begin{equation}
\int_{\{t<\left\vert u\right\vert <t+\kappa\}}\Phi\left(  \nabla
u\right)
dx\geq\int_{\{t<u^{\bigstar}<t+\kappa\}}\Phi_{\blacklozenge}(
\big\vert \nabla u^{\bigstar}\big\vert )  dx \label{A}.
\end{equation}
By (H1), we get
\begin{align}
\frac{1}{\kappa}\int_{\{t<\left\vert u\right\vert \leq
t+\kappa\}}\Phi\left( \nabla u\right)  dx  &
\leq\frac{1}{\kappa}\int_{\{t<\left\vert u\right\vert \leq
t+\kappa\}}(f(x)-g(x,u))\left(  \left\vert u\left(  x\right)
\right\vert
-t\right)  \text{sign}\left(  u\left(  x\right)  \right)  \text{ }%
dx\label{B}\\
&  +\int_{\{\left\vert u\right\vert
>t+\kappa\}}(f(x)-g(x,u))\,\text{sign}\left( u\left(  x\right)
\right)  \text{ }dx\,.\nonumber
\end{align}
We observe that, since $f$ is nonnegative, standard arguments assure that
$u\geq 0$. So, by (H4) and using the monotonicity of $b$, we obtain%
\begin{equation}
\int_{\{|u|>t\}}g(x,u(x))\text{ sign}\;u\text{
}dx\geq\int_{\{|u|>t\}}b(u(x))\text{
sign}\;u\text{ }dx=\int_{0}^{\mu_{u}(t)}b\,(u^{\ast}(s))\;ds. \label{b_1}%
\end{equation}
Letting $\kappa\rightarrow0^{+}$ in (\ref{B}), by (\ref{A}),
(\ref{b_1}) and using the properties of rearrangement, we get
\begin{equation}
0\leq
-\frac{d}{dt}\int_{\{u^{\bigstar}>t\}}\Phi_{\blacklozenge}\left(
\big\vert \nabla u^{\bigstar}\big\vert \right)
dx\leq-\frac{d}{dt}\int_{\{\left\vert u\right\vert >t\}}\Phi\left(
\nabla u\right)  \;dx\leq\int_{0}^{\mu_{u}\left(
t\right)  }(f^{\ast}\left(  s\right)  -b(u^{\ast}(s)))\;ds\,. \label{C}%
\end{equation}
\newline We claim that
\begin{equation}
1\leq\frac{-\mu_{u}^{\prime}\left(  t\right) }{N\omega_{N}^{1/N}
\mu_{u}\left(  t\right)^{1/N^{\prime}}}\Psi_{\blacklozenge}%
^{-1}\left(
\frac{-\frac{d}{dt}\int_{\{u^{\bigstar}>t\}}\Phi_{\blacklozenge
}\left( |\nabla u^{\bigstar}|\right)
dx}{N\omega_{N}^{1/N}\mu_{u}\left(
t\right)  ^{1/N^{^{\prime}}}}\right)\qquad \hbox{for \emph{a.e.}}\;\;t>0\,.  \label{Cianchi_Dis}%
\end{equation}
Indeed, since $u$ and $u^\bigstar$ are equimeasurable, the Jensen
inequality assures that
\begin{equation*}
\Phi_\blacklozenge \left (\frac{\frac 1h \int_{\{ t< u^\bigstar<
t+h\}} |\nabla u^\bigstar| \; dx}{\frac{\mu_{u} (t) - \mu_{u} (t+
h)}{h}}\right ) \leq \frac{\frac 1h \int_{\{ t< u^\bigstar< t+h\}}
\Phi_\blacklozenge (|\nabla u^\bigstar|) \; dx}{\frac{\mu_{u} (t) -
\mu_{u} (t+ h)}{h}}\qquad \hbox{for}\;\;t \;\; \hbox{\rm and}\;\; h
>0\,,
\end{equation*}
that is
\begin{equation*}
\frac{\frac 1h \int_{\{ t< u^\bigstar< t+h\}} |\nabla u^\bigstar| \;
dx}{\frac{\mu_{u} (t) - \mu_{u} (t+ h)}{h}} \leq
\Psi_\blacklozenge^{-1} \left (\frac{\frac 1h \int_{\{ t<
u^\bigstar< t+h\}} \Phi_\blacklozenge (|\nabla u^\bigstar|) \;
dx}{\frac{1}{h}\int_{\{ t< u^\bigstar< t+h\}}|\nabla u^\bigstar|\,
dx }\right)\qquad \hbox{for}\;\;t \;\; \hbox{\rm and}\;\; h
>0\,.
\end{equation*}
Finally, coarea formula yields \eqref{Cianchi_Dis}.
\\
 Combining  (\ref{C}) and \eqref{Cianchi_Dis}, we obtain
\[
1\leq\frac{-\mu_{u}^{\prime}\left(  t\right) }{N\omega_{N}^{1/N}
\mu_{u}\left(  t\right)^{1/N^{\prime}}}\Psi_{\blacklozenge}%
^{-1}\left(  \displaystyle\frac{\int_{0}^{\mu_{u}\left(  t\right)
}(f^{\ast }\left(  s\right)
-b(u^{\ast}(s)))\;ds}{N\omega_{N}^{1/N}\mu _{u}\left(
t\right)^{1/N^{\prime}}}\right)  \qquad \hbox{for
\emph{a.e.}}\;\;t>0\,.
\]
By standard arguments, it follows that%

\begin{equation}
-(u^{\ast}(s))^{^{\prime}}\leq\frac{1}{N\omega_{N}^{1/N}\,s^{1/N^{\prime}}%
}\Psi_{\blacklozenge}^{-1}\left(  \displaystyle\frac{\mathcal{F}%
(s)-\mathcal{B}(s)}{N\omega_{N}^{1/N}s^{1/N^{\prime}}}\right)
\qquad \hbox{for }\;s\in (0, |\Omega|)\,. \label{E:10}%
\end{equation}
By \eqref{3:11}, the derivative of $\mathcal{B}$ equals
\begin{equation}
\mathcal{B}^{\prime}(s)=b(u^{\ast}(s))\qquad\hbox{for
\emph{a.e.}}\;s\in(0,|\Omega|). \label{3:13}%
\end{equation}
Equations \eqref{3:11}, \eqref{E:10} and \eqref{3:13} yield
\begin{equation}
\left\{
\begin{array}
[c]{lll}%
N\omega_{N}^{1/N}\;s^{1/N^{\prime}}\Psi_{\blacklozenge}\left(
N\omega _{N}^{1/N}\;s^{1/N^{\prime}}\left(
-\displaystyle\frac{d}{ds}\left( \gamma\left(
\mathcal{B}^{\prime}(s)\right)  \right)  \right)  \right)
+\mathcal{B}(s)\leq\mathcal{F}(s) & \qquad\hbox{for
$s\in(0, |\Omega|])$} & \\
&  & \\
\mathcal{B}(0)=0,\;\;\mathcal{B}^{\prime}(|\Omega|)=0\,, &  &
\end{array}
\right.  \label{E:10'}%
\end{equation}
where $\gamma(\sigma)=b^{-1}(\sigma),$ the inverse function of $b$.

Now let us consider problem \eqref{Prob_sym_zero}. The solution $v$
to problem \eqref{Prob_sym_zero} is unique and the symmetry of data
assures that $v(x)=v(|x|)$, \emph{i.e.} $v$ is positive and radially
symmetric. Moreover, setting $s=$ $\omega_{N}\left\vert x\right\vert
^{N}$ and $\widetilde {v}\left(  s\right)  =v\left(  \left(
s/\omega_{N}\right)  ^{1/N}\right)  ,$ we get, for all $s\in\left[
0,\left\vert \Omega\right\vert \right]  $,
\begin{equation}\label{=}
-\,\frac{\Phi_{\blacklozenge}\left( | \widetilde{v}^{\prime}\left(
s\right)| N\omega_{N}^{1/N}\;s^{1/N^{\prime}}\right)}{
\widetilde{v}^{\prime }(s)}=\int_{0}^{s}\left(
\widetilde{f}^{\ast}\left(  \sigma\right)
-b(\widetilde{v}(\sigma))\right)  d\sigma.
\end{equation}

\noindent Arguing as in \eqref{C}, we have that the integral on the
right-hand side of \eqref{=} is
positive and then this fact guarantees that%

\[
\,\frac{\Phi_{\blacklozenge}\left( | \widetilde{v}^{\prime}\left(
s\right)| N\omega_{N}^{1/N}\;s^{1/N^{\prime}}\right)}{
\widetilde{v}^{\prime }(s)}\leq0.
\]
This means  that $v(x)=v^{\bigstar}(x)$, since
$\Phi_{\blacklozenge}$ a positive function.

By the properties of $v$, we can repeat arguments used to prove
\eqref{E:10} replacing all the inequalities by equalities (see also
\cite[Lemma 1.32]{Dz85}) and obtaining
\[
-(v^{\ast}(s))^{^{\prime}}=\frac{1}{N\omega_{N}^{1/N}\,s^{1/N^{\prime}}}%
\Psi_{\blacklozenge}^{-1}\left(  \displaystyle\frac{\widetilde{\mathcal{F}%
}(s)-\widetilde{\mathcal{B}}(s)}{N\omega_{N}^{1/N}s^{1/N^{\prime}}}\right),
\]
and thus
\begin{equation}
\left\{
\begin{array}
[c]{lll}%
N\omega_{N}^{1/N}\;s^{1/N^{\prime}}\Psi_{\blacklozenge}\left[
N\omega _{N}^{1/N}\;s^{1/N^{\prime}}\left(
-\displaystyle\frac{d}{ds}\left( \gamma\left(
\widetilde{\mathcal{B}}^{\prime}(s)\right)  \right)  \right) \right]
+\widetilde{\mathcal{B}}(s)=\widetilde{\mathcal{F}}(s) &
\qquad\hbox{for
$s\in(0, |\Omega|])$} & \\
&  & \\
\widetilde{\mathcal{B}}(0)=0,\;\;\widetilde{\mathcal{B}}^{\prime}%
(|\Omega|)=0\,. &  &
\end{array}
\right.  \label{E:10''}%
\end{equation}
Since
$\mathcal{B},\widetilde{\mathcal{B}}\in\mathcal{C}([0,|\Omega|])$,
there exists $\overline{s}\in(0,|\Omega|)$ such that
\begin{equation*}
\Vert(\mathcal{B}-\widetilde{\mathcal{B}})_{+}\Vert_{L^{\infty}(0,|\Omega
|)}=(\mathcal{B}-\widetilde{\mathcal{B}})(\overline{s}). \label{E:30}%
\end{equation*}
We argue by contradiction. Suppose that
\begin{equation}
(\mathcal{B}-\widetilde{\mathcal{B}})(\overline{s})>\Vert(\mathcal{F}-\widetilde
{\mathcal{F}})_{+}\Vert_{L^{\infty}(0,|\Omega|)}\,. \label{E:31}%
\end{equation}
If $\overline{s}<|\Omega|$, combining (\ref{E:10'}) and
(\ref{E:10''}) yields
\begin{align}
&  N\omega_{N}^{1/N}\;s^{1/N^{\prime}}\left[
\Psi_{\blacklozenge}\left( N\omega_{N}^{1/N}\;s^{1/N^{\prime}}\left(
-\displaystyle\frac{d}{ds}\left( \gamma\left(
\mathcal{B}^{\prime}(s)\right)  \right)  \right)  \right)
-\Psi_{\blacklozenge}\left(
N\omega_{N}^{1/N}\;s^{1/N^{\prime}}\left(
-\displaystyle\frac{d}{ds}\left(  \gamma\left(  \widetilde{\mathcal{B}%
}^{\prime}(s)\right)  \right)  \right)  \right)  \right] \label{E:35}\\
&  \leq\mathcal{F}(s)-\widetilde{\mathcal{F}}(s)+\widetilde{\mathcal{B}%
}(s)-\mathcal{B}(s).\nonumber
\end{align}
By (\ref{E:31}), it follows that
\begin{equation}
\mathcal{F}(s)-\widetilde{\mathcal{F}}(s)+\widetilde{\mathcal{B}%
}(s)-\mathcal{B}(s)\leq\Vert(\mathcal{F}-\widetilde{\mathcal{F}})_{+}%
\Vert_{L^{\infty}(0,|\Omega|)}-(\mathcal{B}-\widetilde{\mathcal{B}%
})(s)<0\qquad\hbox{for}\;s\in(\overline{s}-\varepsilon,\overline{s}+\varepsilon).
\label{E:36}%
\end{equation}
We set
\begin{equation*}
Z=\mathcal{B}-\widetilde{\mathcal{B}}\in
W^{2,\infty}(\overline{s}-\varepsilon
,\overline{s}+\varepsilon). \label{E:32}%
\end{equation*}
Then,
\begin{equation}
\gamma\left(  \mathcal{B}^{\prime}(s)\right)  -\gamma\left(
\widetilde
{\mathcal{B}}^{\prime}(s)\right)  =Z^{\prime}(s)\,c(s), \label{E:33}%
\end{equation}
where
\begin{equation*}
c(s)=\int_{0}^{1}\gamma^{\prime}\left(  \tau\mathcal{B}^{\prime}%
(s)+(1-\tau)\widetilde{\mathcal{B}}^{\prime}(s)\right)  \;d\tau>0.
\label{E:34}%
\end{equation*}
As a consequence of (\ref{E:35}) and (\ref{E:36}), we obtain
\begin{equation}
\Psi_{\blacklozenge}\left( N\omega_{N}^{1/N}\;s^{1/N^{\prime}}\left(
-\displaystyle\frac{d}{ds}\left(  \gamma\left(  \mathcal{B}^{\prime
}(s)\right)  \right)  \right)  \right)  -\Psi_{\blacklozenge}\left(
N\omega_{N}^{1/N}\;s^{1/N^{\prime}}\left(
-\displaystyle\frac{d}{ds}\left( \gamma\left(
\widetilde{\mathcal{B}}^{\prime}(s)\right)  \right)  \right)
\right)  <0. \label{E:37}%
\end{equation}
By the monotonicity of $\Psi_{\blacklozenge}^{-1}$ and by
(\ref{E:37}), we get
\[
-\displaystyle\frac{d}{ds}\left(  \gamma\left(  \mathcal{B}^{\prime
}(s)\right)  \right)  <-\displaystyle\frac{d}{ds}\left( \gamma\left(
\widetilde{\mathcal{B}}^{\prime}(s)\right)  \right).
\]
Using (\ref{E:33}), we can conclude that
\begin{equation}
-\frac{d}{ds}\left(  c(s)Z^{\prime}(s)\right)  <0\qquad\hbox{for
}\;s\in(\overline{s}-\varepsilon,\overline{s}+\varepsilon), \label{ZZ}%
\end{equation}
which is in contradiction with the assumption that $Z$ has a maximum
in $\overline{s}.$

If $\overline{s}=\left\vert \Omega\right\vert ,$ (\ref{ZZ})\ holds
in $(\left\vert \Omega\right\vert -\varepsilon,\left\vert
\Omega\right\vert )$ and then $Z^{\prime}(\left\vert
\Omega\right\vert )>0$. But we know that $Z^{\prime }(\left\vert
\Omega\right\vert )=0.$
%

\subsection{Generalization }
We point out that some generalizations of Theorem \ref{teomain}
hold. Let $b_{1},b_{2}$ be two continuous strictly increasing
functions. We say that $b_{1}$ is \textit{weaker} than
$b_{2}$ and we write%
\begin{equation*}
b_{1}\prec b_{2}, \label{b minore b tilde}%
\end{equation*}
if they have the same domains and there exists a
contraction\footnote{By contraction, we mean $\left\vert
\rho(a)-\rho(b)\right\vert \leq\left\vert a-b\right\vert$ for $a, b
\in \R$.} $\rho: \mathbb{R}
\rightarrow%
\mathbb{R}$ such that $b_{1}=\rho\circ b_{2}.$

It is possible to estimate the difference between the concentration
of a solution $u\in V_{0}^{1,\Phi}(\Omega)$ to problem \eqref{Prob}
and the
concentration of the solution $v\in V_{0}^{1,\Phi_{\blacklozenge}}%
(\Omega^{\bigstar})$ to the following problem
\begin{equation*}
\left\{
\begin{array}
[c]{ll}%
-\operatorname{div}\left(  \!\dfrac{\Phi_{\blacklozenge}\left(
\left\vert \nabla v\right\vert \right)  }{\left\vert \nabla
v\right\vert ^{2}}\nabla v\!\right)
+\widetilde{b}(v)=\widetilde{f}(x) &
\hbox{in $\Omega^\bigstar$}\\
& \\
v=0 & \hbox{on $\partial\Omega^\bigstar$},
\end{array}
\right.  \label{prob_sym_zero_2}%
\end{equation*}
under the assumptions
(H$^\bigstar$1),(H$^\bigstar$3)-(H$^\bigstar$5) and with
$\widetilde{b}$ a continuous strictly increasing function such that
$\widetilde {b}(0)=0$ and $\widetilde{b}^{-1}\prec b^{-1}$.

More precisely, we get
\[
\Vert(\mathcal{B}-\widetilde{\mathcal{B}})_{+}\Vert_{L^{\infty}(0,|\Omega
|)}\leq\Vert(\mathcal{F}-\widetilde{\mathcal{F}})_{+}\Vert_{L^{\infty
}(0,|\Omega|)},
\]
where $\mathcal{F}$ and $\widetilde{\mathcal{F}}$ are defined as in
(\ref{3:11'}) and, with abuse of notation, $\mathcal{B}$ and
$\widetilde{\mathcal{B}}$ are defined by
\[
\mathcal{B}(s)=\int_{0}^{s}b(u^{\ast}(t))\;dt\qquad\widetilde{\mathcal{B}%
}(s)=\int_{0}^{s}\widetilde{b}(v^{\ast}(t))\;dt
\]
for $s\in(0,|\Omega|]$. In particular, if we suppose that
\[
\mathcal{F}(s)\leq\widetilde{\mathcal{F}}(s)\qquad\hbox{for
any}\,\;s\in\lbrack0,|\Omega|],
\]
then%
\[
\mathcal{B}(s)\leq\widetilde{\mathcal{B}}(s)\qquad\hbox{for
any}\,\;s\in\lbrack0,|\Omega|].
\]

Finally, we stress that in Theorem \ref{teomain} we can replace
assumption that
 $b$ is a continuous and non decreasing function such that $b(0)=0$
by the following more general condition that $b$  is a maximal
monotone graphs in $ \mathbb{R}^{2}$ such that $b(0)\ni0.$ Indeed, a
maximal monotone graph is a natural generalization of the concept of
monotone non-decreasing real function; moreover, the inverse of a
maximal monotone graph (that appears in the proof of Theorem
\ref{teomain}) is again a maximal monotone graph (see \cite{V} for
more details). Results in this order of idea are contained in
\cite{AdBF3} when $\Phi$ is given by \eqref{Fi pi}.

\section{Appendix: existence and uniqueness for problem \eqref{Prob_sym_zero}}

In this Appendix, we focus our attention on  the existence and
uniqueness of a solution to the symmetrized problem
\eqref{Prob_sym_zero}.
\\
Let us consider boundary value problems for elliptic equation of the
form
\begin{equation}
\left\{
\begin{array}
[c]{ll}%
-\operatorname{div}\left( \textbf{A}\left( x,  \nabla v
\right)\right)
+b(v)= h(x) & \hbox{in $\Omega$}\\
& \\
v=0 & \hbox{on $\partial\Omega$}\,.
\end{array}
\right.  \label{Prob_GM}%
\end{equation}
It is well-known that Br\'{e}zis and Browder in \cite{B}  studied
the existence for this type of problems when the coefficients of
$\textbf{A} =(A_1, \ldots, A_m)$ have polynomial growth and $b$
satisfies a sign condition.
\\
Here,  we suppose that the principal operator in \eqref{Prob_GM}
have general growth not necessarily  of polynomial type, but
governed by a Young function that does not necessarily verify the
$\triangle_2$-condition. Then, this kind of problems are treated in
the framework of the Orlicz-Sobolev spaces  and the difficulty comes
from the fact that, in general, these spaces are not reflexive.
\\
Questions related to the existence of the solution to such a problem
in the form \eqref{Prob_GM} have been studied, for example, in
\cite{GM}.
For the convenience of the reader, we recall the following
definition and the existence result.

\medskip
An open subset $\Omega$ of $\rN$ has the segment property if there
exist a locally finite open covering $\{O_i\}$ of $\partial \Omega$
and corresponding vectors $\{\eta_i\}$ such that $ x+\tau\eta_i \in
\Omega$ for every $x\in \bar{\Omega}\cap O_i$ and $0<\tau<1$,
\emph{i.e.} its boundary is locally the graph of a continuous
function.
\begin{theorem}\label{teo_existence}
Let $\Omega$ be a bounded open subset of $\rN$  having the segment
property. Let us consider problem \eqref{Prob_GM} such that the
following conditions hold:
\begin{itemize}
\item[$(A_1)$]{Every component $A_i$ is a real-valued function
defined in $\Omega\times\R^N$  which is measurable in $x$ for fixed
$\xi$ and continuous in $\xi$ for a.e. fixed $x$.
\item[$(A_2)$] There exist a one-dimensional Young function $M$,
a function $\nu\in E_{M_{\bullet}}(\Omega)$ and two constants $c_1$
and $c_2$ such that
$$|A_i(x,\xi)|\leq \nu(x)+c_1 {M_{\bullet}}^{-1}M(c_2|\xi|)$$
for a.e. $x \in \Omega$ and all $\xi \in \R^N$.}
\item[$(A_3)$]{
$\displaystyle\left (\textbf{A}(x,\xi_1) - \textbf{A}(x,\xi_2)
\right)\cdot \left(\xi_1 - \xi_2\right) > 0\qquad \hbox{\rm
\emph{for a.e.} $x \in \Omega$ \emph{ and  every} $\xi_1 \neq \xi_2$
\emph{in} $\rN$}\,.$}
\item[$(A_4)$] {There exist two functions $\textbf{b}\in
 (E_{M_{\bullet}}(\Omega))^N$, $b_0\in L^1(\Omega)$ and
 two constants $d_1$ and $d_2$ such that for some fixed element
 $\varphi \in W_0^1E_M(\Omega)\cap L^{\infty}(\Omega)$
$$\textbf{A}(x,\xi)\cdot (\xi-\nabla \varphi)\geq d_1
M(d_2|\xi|)-\textbf{b}\cdot \xi-b_0(x)$$ for a.e. $x \in \Omega$ and
all $\xi \in \R^N$.}
\end{itemize}
Let $b: \R \va \R$ be a Carath\'{e}odory function such that $b(v)\,
v\geq 0$ for all $v$. Then, given $h\in E_{M_\bullet}(\Omega)$,
there exists a weak solution $v\in W_0^1 L_M(\Omega)\cap\left\{v \in
W_0^1 L_M(\Omega): |\textbf{A}(x,\nabla v)|\in
L_{M_{\bullet}}(\Omega)\right\}$ to problem \eqref{Prob_GM}.

\end{theorem}
\bigskip

We stress that, when the components of $\textbf{A}$ have a
polynomial growth, \textit{i.e.} when $M(t)=|t|^p$ for $p>1$, then
conditions (A1)-(A4) reduce to the classical conditions of the
Leray-Lions operators.

The previous existence result runs for $\textbf{A}(x, \nabla v)=
\Big(\frac{m(|\nabla v|)}{|\nabla v|}\nabla v\Big)$ and
$M(t)=\int_0^{|t|} m(s) \,ds$, where $m:[0,+\infty)\rightarrow [0,+
\infty)$ is an increasing, right continuous function, $m(t)=0$ if
and only if $t=0$ and $m(t)\rightarrow +\infty$ as $t \rightarrow
+\infty$.

Then, it can be applied to the symmetrized problem
\eqref{Prob_sym_zero} with $\Omega=\Omega^\bigstar$,
$m(t)=\frac{\Phi_\blacklozenge (t)}{t}$ and $M(t)=\Theta (t)$, where
$\Theta$ is defined as in \eqref{theta}. In fact, the assumptions
$(A_1)-(A_4)$ are fulfilled with respect to some fixed $\varphi\in
W_0^1E_{\Theta}(\Omega)\cap L^{\infty}(\Omega)$. More precisely,
condition $(A_1)$ is obvious since \eqref{psi incrising} holds.
Condition $(A_2)$ follows  by \eqref{theta} and \eqref{dis CB}
applied to the Young function $\Theta$:
\begin{align} \nonumber
\Theta(2t)\geq \int_t^{2t} \frac{\Phi_\blacklozenge (\tau)}{\tau}\;
d\tau \geq
\Phi_\blacklozenge (t) \geq
\Theta_{\bullet}\left(\frac{\Phi_\blacklozenge (t)}{t}\right) \qquad
\hbox{\rm for}\;\; t\geq 0,
\end{align}
\textit{i.e.}
\begin{equation*}
\frac{\Phi_\blacklozenge (t)}{t} \leq \left(
\Theta_{\bullet}\right)^{-1}\left(\Theta(2t)\right)\qquad \hbox{\rm
for}\;\; t\geq 0.
\end{equation*}
\\
Condition $(A_3)$ is a consequence of the strictly convexity
\footnote{The function $\Theta(|\xi|)$ verifies the strictly convex
condition as $N-$variable function. This follows combining the
strictly convexity and the strictly monotonicity (as function of one
variable) in $[0,+\infty)$ of $\Theta$ with the triangular
inequality of modulus.
 We have only to check that $\Theta$ is strictly convex as one variable function.
  Indeed, the strictly increasing of  $\Phi_\blacklozenge$ implies that $s_0 =0$
   in \eqref{s0} and then $\frac{\Phi_\blacklozenge(s)}{s}$ is strictly increasing in $[0,+\infty)$. }
 of the
function $\xi\in \rN\rightarrow\Theta(|\xi|)\in\R$. Indeed, for all
$\xi_1, \xi_2\in \rN$ with $\xi_1\neq \xi_2$, we get
\begin{equation*}
\begin{aligned}\displaystyle\left(\frac{\Phi_\blacklozenge (|\xi_1|)}{|\xi_1|^2}
\xi_1 - \frac{\Phi_\blacklozenge (|\xi_2|)}{|\xi_2|^2} \xi_2
\right)\cdot \left(\xi_1 - \xi_2\right)
=\displaystyle\left(\nabla\Theta(|\xi_1|)-
\nabla\Theta(|\xi_2|)\right)\cdot \left(\xi_1 -
\xi_2\right)\\>-\left[\Theta(|\xi_2|)-\Theta(|\xi_1|)\right]-
\left[\Theta(|\xi_1|)-\Theta(|\xi_2|)\right]=0.
\end{aligned}\end{equation*}
By applying Proposition \ref{equalityYoung} \footnote{For all
$\xi,\eta \in \rN$, with $|\eta|=1$, we have that
$\frac{\partial}{\partial\eta}\Theta(|\xi|)=\Theta'(|\xi|)\nabla\xi\cdot\eta$
is linear in $\eta$, where $\Theta'(t)$ is the derivative of
$\Theta(t)$.
 }  to the Young function $\Theta$ and the Young inequality \eqref{Young},
it follows that
\begin{align} \nonumber
&\left (\frac{\Phi_\blacklozenge (|\xi|)}{|\xi|^2} \xi\right) \left
(\xi - \nabla \varphi (x)\right) \geq \Theta_{\bullet}\left
(\frac{\Phi_\blacklozenge (|\xi|)}{|\xi|} \right) + \Theta (|\xi|) -
\Theta_{\bullet}\left (\frac{\Phi_\blacklozenge (|\xi|)}{|\xi|}
\right) -\Theta (|\nabla \varphi (x)|) \\ \nonumber \cr &= \Theta
(|\xi|) - \Theta (|\nabla \varphi (x)|), \nonumber
\end{align}
namely condition $(A_4)$ with constants $d_1=d_2=1$,
$\textbf{b}(x)=\textbf{0} $ and $b_0(x)= \Theta (|\nabla \varphi
(x)|)$ in $L^1(\Omega)$.
\\

%
%
%
%
%
%
%
%


Finally, we verify the uniqueness of the solution 
to the symmetrized problem \eqref{Prob_sym_zero}. It follows in a
standard way using the monotonicity of zero order term. Indeed, let
$v_1$ and $v_2$ be two different solutions to \eqref{Prob_sym_zero}
and let us choose
$T_{\delta}\left(\left(v_1-v_2\right)_+\right)/\delta$
  as test function in the difference of equations,
  where $T_{\delta}(t):= \min \{\delta, \max \{t, -\delta\}\}$.
  Then,
  \begin{align}\label{4.10}
  \int_{\Omega^{\bigstar}} \Big[\Phi_\blacklozenge (|\nabla v_1|)\frac{\nabla v_1}{|\nabla v_1|^2}
  -\Phi_\blacklozenge (|\nabla v_2|)\frac{\nabla v_2}{|\nabla v_2|^2} \Big]&
  \frac{\nabla T_\delta \left(\left(v_1-v_2\right)_{+}\right)}
  {\delta} \; dx \\
  &+\int_{\Omega^{\bigstar}} \Big[b(v_1)-b(v_2)\Big] \frac{\nabla T_\delta \left(\left(v_1-v_2\right)_{+}\right)}
  {\delta} \; dx=0 . \nonumber
  \end{align}
By (A3), \eqref{4.10} and the definition of
$\left(v_1-v_2\right)^+$, we obtain
\begin{equation*}\label{4.12}
\int_{\Omega^{\bigstar}} [b(v_1) - b(v_2)] \frac{\nabla T_\delta
\left(\left(v_1-v_2\right)_{+}\right)}
  {\delta} \; dx \leq 0,
  \end{equation*}
 and, passing to the limit with $\delta \va 0$
and using the monotonicity of $b$, it follows $ |\{v_1-v_2>0\}| = 0.
$ Replacing $v_2$ by $v_1$, the assert follows.

Similar arguments show that the solution $v$ to problem
\eqref{Prob_sym_zero} is nonnegative since $f$ is nonnegative.

\paragraph*{Acknowledgements}

This work has been partially supported by GNAMPA of the Italian
INdAM (National Institute of High Mathematics) and ``Programma
triennale della Ricerca dell'Universit\`{a} degli Studi di Napoli
``Parthenope'' - Sostegno alla ricerca individuale 2015-2017''.

\end{document}